\documentclass{article}
\usepackage{a4,amsmath,amssymb,amsfonts,amsthm}
\usepackage[all]{xy}
\usepackage[latin1]{inputenc}

\begin{document}
\bibliographystyle{plain}

%
%
 
%
\def\M#1{\mathbb#1}     
\def\B#1{\bold#1}       
\def\C#1{\mathcal#1}    
\def\E#1{\scr#1}        
\def\mR{\M{R}}           %
\def\mZ{\M{Z}}           %
\def\mN{\M{N}}           
\def\mQ{\M{Q}}           %
\def\mC{\M{C}}    		%
\def\mP{\M{P}}

%
%
\gdef\beginProof{\par{\bf Proof. }}
\gdef\endProof{$\Box$\par}  
\gdef\ari#1{\widehat{#1}}
\gdef\pr{^{\prime}}
\gdef\prpr{^{\prime\prime}}
\gdef\mtr#1{\overline{#1}}
\gdef\ra{\rightarrow}
\gdef\Bbb{\bf }
\gdef\P1{{\Bbb P}^{1}_{D}}
\gdef\dbd{{{i\over 2\pi}\partial\overline{\partial}}}
\gdef\a{\alpha}
\gdef\ca{\ch(\alpha)}
\gdef\ttmk#1{\widetilde{\theta}^{k}(#1)}
\gdef\tdm#1{\theta^{k}({#1}^{\vee})^{-1}}
\gdef\td#1{\theta^{k}({#1}^{\vee})}
\gdef\rl{{\Lambda}}
\gdef\CT{CT}
\gdef\refeq#1{(\ref{#1})}
\gdef\div{{\rm div}}


\gdef\CH{{\rm CH}}
\gdef\Spec{{\rm Spec}}
\gdef\ch{{\rm ch}}
\gdef\Im{{\rm Im}}
\gdef\Td{{\rm Td}}
\gdef\td{{\rm td}}
\gdef\red{{\rm red}}

 \newtheorem{theor}{Theorem}[section]
 \newtheorem{prop}[theor]{Proposition}
 \newtheorem{cor}[theor]{Corollary}
 \newtheorem{lemma}[theor]{Lemma}
 \newtheorem{sublem}[theor]{sublemma}
 \newtheorem{defin}[theor]{Definition}
 \newtheorem{conj}[theor]{Conjecture}

 \author{Henri Gillet\footnote{Department of Mathematics, 
 University of Illinois at Chicago, 
 Box 4348, Chicago IL 60680, USA} , 
 Damian Rössler\footnote{Institut de Math\'ematiques de Jussieu, CNRS,
 Case Postale 7012,
 2 place Jussieu,
 F-75251 Paris Cedex 05, France}\ \ and C. Soulé\footnote{CNRS and IHES, 
  35 Route de Chartres, F-91440 
 Bures-Sur-Yvette, France}}
 \title{An arithmetic Riemann-Roch theorem in higher degrees}
 \maketitle
 \date
 \abstract{We prove an analogue in Arakelov geometry of the Grothendieck-Riemann-Roch theorem.}
\begin{flushleft}
\end{flushleft}
\begin{flushleft}
{\bf Mathematics Subjects Classification 2000}: 14G40, 14C40, 58J52
\end{flushleft}

 \parindent=0pt

\section{Introduction}

Recall that the Grothendieck-Riemann-Roch 
theorem (see for instance \cite[par. 20.1]{Fulton}) says that, if $Y$ and $B$ are regular schemes which are quasi-projective and flat over the spectrum  $S$ of a Dedekind domain and $g:Y\ra B$ is a flat and 
projective 
$S$-morphism, then the diagram
$$
\xymatrix{
K_0(Y)\ar[rr]^{\Td(g)\cdot\ch}\ar[d]_{g_*} & & \CH^\cdot(Y)_\mQ\ar[d]^{g_*\ \ \ \ \ \ \ \ \ \ {\rm (GRR)}}\\
K_0(B)\ar[rr]^{\ch} & & \CH^\cdot(B)_\mQ
}
\label{diaggrr}
$$
commutes. Here $K_0(Y)$ (resp. $K_0(B)$) is the Grothendieck group of locally 
free sheaves on $Y$ (resp. on $B$). The group $\CH^\cdot(Y)$ 
(resp. $\CH^\cdot(B)$) is the Chow group of cycles modulo rational equivalence 
on $Y$ (resp. $B$). The symbol $g_*$ refers to the push-forward map 
in the corresponding theory. 
The symbol $\ch$ refers to the Chern 
character, which on each regular and quasi-projective $S$-scheme 
is a ring morphism from the Grothendieck group to the Chow group 
tensored with $\mQ$. The element $\Td(g)$ is the Todd class 
of the virtual relative tangent bundle. In words, the Grothendieck-Riemann-Roch theorem 
implies that the Chern character does not commute with the push-forward maps but 
that this commutation can be obtained after multiplication of the Chern 
character with the Todd class of the virtual relative tangent bundle. 

All the objects mentionned in the previous paragraph have extensions  
to Arakelov theory. Arakelov theory is an extension of scheme-theoretic 
geometry over the integers, where everything in sight is equipped 
with an analytic datum on the complex points of the scheme. 
This means that
 there will be "forgetful" maps from the Arakelov-theoretic 
 Grothendieck groups, Chow groups etc. to the corresponding 
 classical objects. We refer to \cite{Lang3} and \cite{Soule}
 for an introduction to this subject, which originated in Arakelov's paper 
 \cite{Ara} and was later further developped by several people. 
 
Now keep the same hypotheses as in the diagram (GRR) but 
suppose in addition that $S=\Spec\ \mZ$ and that $g$ is smooth over $\mQ$. 
We shall show that 
there exists a commutative diagram 
$$
\xymatrix{
\ari{K}_0(Y)\ar[rrr]^{\ari{\Td}(g)\cdot(1-R(Tg_\mC))\cdot\ari{\ch}}\ar[d]_{g_*} & & &\ari{\CH}^\cdot(Y)_\mQ\ar[d]^{g_*\ \ \ \ \ \ \ \ \ \ {\rm (ARR)}}\\
\ari{K}_0(B)\ar[rrr]^{\ari{\ch}} & & &\ari{\CH}^\cdot(B)_\mQ
\label{diagarr}
}
$$
where the objects with hats ($\ \ari{\cdot}\ $) are the extensions of the 
corresponding objects to Arakelov theory. The class 
$1-R(Tg_\mC)$ is an exotic cohomology class which has no classical 
analog. This diagram fits in a  three-dimensional commutative 
diagram
$$
\xymatrix{
K_0(Y)\ar[rr]^{\Td(g)\cdot\ch}\ar[dd]_{g_*} & &\CH^\cdot(Y)_\mQ\ar@{..>}'[d][dd]_{\ g_*}&\\
&\ari{K}_0(Y)\ar[dd]_>>>>>>>{g_*}\ar@{->>}[ul]\ar[rr]^{\ari{\Td}(g)(1-R(Tg_\mC))\cdot\ari{\ch}}& 
&\ari{\CH}^\cdot(Y)_\mQ\ar[dd]_{g_*}\ar@{->>}[ul]\\          
K_0(B)\ar'[r][rr]^>>>>>>>>>{\ch} & &\CH^\cdot(B)_\mQ&\\
&\ari{K}_0(B)\ar@{->>}[ul]\ar[rr]^{\ari{\ch}} & &\ari{\CH}^\cdot(B)_\mQ\ar@{->>}[ul]
}
$$
where the various forgetful arrows $\twoheadrightarrow$ are surjective and their 
kernels 
are spaces of differential forms 
on the complex points of the corresponding schemes. 
The assertion that the diagram (ARR) commutes shall 
henceforth be referred to as the {\it arithmetic Riemann-Roch theorem}. 
The precise statement is given in Theorem \ref{arr} below.

Our proof of the arithmetic Riemann-Roch theorem 
combines the classical technique of 
proof of the Grothendieck-Riemann-Roch theorem with deep results of Bismut 
and his coworkers in local index theory. 

The history of the previous work on this theorem and its variants is as follows. 
In \cite{Faltings2} Faltings proved a variant of the theorem for surfaces. 
In \cite[Th. 7, ii)]{GS8}, Gillet and Soulé proved a degree one version of the theorem 
(see after Theorem \ref{arr} for a precise statement). 
In his book \cite{Faltings}, Faltings outlined an approach to the proof of 
Theorem \ref{arr}, which is not based on Bismut's work. 
In \cite[par. 8]{R2},  Rössler proved a variant of Theorem \ref{arr}, where the Chow groups 
are replaced by graded $\ari{K}_0$-groups (in the spirit of \cite{SGA6}). This variant is 
a formal consequence of the Riemann-Roch theorem for the Adams operations 
acting on arithmetic ${\rm K}_0$-groups, which is the main result of 
\cite{R2}.  Finally,
in his unpublished thesis \cite{Zha}, Zha obtained a general arithmetic Riemann-Roch theorem 
which does not involve analytic torsion.

A variant of Theorem \ref{arr} was conjectured in \cite[Conjecture 3.3]{GS5}. 
That conjecture is a variant of Theorem \ref{arr} in the sense that a definition of 
the push-forward map is used there which may differ from 
the one used here. A precise comparison has yet to be made. 
In the present setting  
all the morphisms are local complete  intersections, because all the schemes 
are assumed to be regular; 
it is an open problem is to allow singularities at finite places and/or to allow 
more general morphisms. This problem is solved in degree $1$ in \cite[Th. 7, i)]{GS8}. Note that arithmetic Chow theory (like ordinary Chow theory) is not defined outside the category of regular schemes and to tackle the problem of a 
Riemann-Roch theorem for singular schemes, one has to first extend that theory.

The structure of the article is the following. 
In the second section, we recall the definitions of the 
Grothendieck and Chow groups in Arakelov theory and some 
of their basic properties. 
In the third section, we formulate the arithmetic Riemann-Roch theorem. 
In the fourth section, we give a proof of 
the latter theorem. See the beginning of that section for 
a description of the structure of the proof. 

{\bf Acknowledgments.} The second author thanks J.-I. Burgos and K. Köhler 
for interesting discussions related to the contents of this article. 

\section{Arithmetic Grothendieck and Chow groups}

In this section, we shall define extensions of 
the classical Grothendieck and Chow groups to the  framework of Arakelov theory. 

Let 
$X$ be a regular scheme, which is quasi-projective and flat  
over $\mZ$. We shall call such a scheme an {\it arithmetic variety} 
(this definition is more restrictive than the definition given in \cite[sec. 3.2]{GS2}). 
Complex 
conjugation induces an antiholomorphic automorphism $F_{\infty}$ on the 
manifold of complex points $X(\mC)$ of $X$. We shall write 
$A^{p,p}(X)$ for the set of real differential forms $\omega$
of type $p,p$ on $X(\mC)$, which satisfy the 
equation $F_{\infty}^{*}\omega=(-1)^{p}\omega$ and 
we shall write $Z^{p,p}(X)\subseteq A^{p,p}(X)$ for the kernel of the operation 
$d=\partial+\mtr{\partial}$. We also define $\widetilde{A}(X):=\bigoplus_{p\geq 0}(A^{p,p}(X)/(\Im\ \partial+\Im\ \mtr{\partial}))$ and $Z(X):=\bigoplus_{p\geq 0}Z^{p,p}(X)$. A {\it hermitian bundle} $\mtr{E}=(E,h^E)$ is a vector bundle $E$ on $X$, endowed with a 
hermitian metric $h^E$, which is invariant under $F_{\infty}$, 
on the holomorphic bundle $E_\mC$ on $X(\mC)$ 
associated to $E$. We denote by 
$\ch(\mtr{E})$ (resp. $\Td(\mtr{E})$) the representative of the Chern character 
(resp. Todd class) of $E_\mC$ 
associated by the formulae of Chern-Weil to 
the hermitian connection of type $(1,0)$ defined by $h^E$.
Let 
$$
{\cal E}:0\ra E\pr\ra E\ra E\prpr\ra 0
$$ 
be an 
exact sequence of vector bundles on $X$. We shall write $\mtr{\cal E}$ 
for the sequence $\cal E$ and hermitian metrics on $E_\mC\pr$, $E_\mC$ 
and $E_\mC\prpr$ 
(invariant under $F_{\infty}$). To $\mtr{\cal E}$ is 
 associated a secondary, or Bott-Chern 
class $\widetilde{\ch}(\mtr{\cal E})\in\widetilde{A}(X)$ (resp. 
$\widetilde{\Td}(\mtr{\cal E})$). This 
secondary class satisfies the 
equation
$$
\dbd(\widetilde{\ch}(\mtr{\cal E}))=\ch(\mtr{E}'\oplus\mtr{E}'')-\ch(\mtr{E})
$$
(resp.
$$
\dbd(\widetilde{\Td}(\mtr{\cal E}))=\Td(\mtr{E}'\oplus\mtr{E}'')-\Td(\mtr{E}){\ \rm ).}
$$
Here we write $\mtr{E}'\oplus\mtr{E}''$ for the hermitian bundle 
$(E'\oplus E'',h'\oplus h'')$, which is the orthogonal direct sum of 
the hermitian bundles $\mtr{E}'$ and $\mtr{E}''$. 
For the definition of the secondary classes, we refer to \cite[Par. f)]{BGS1}. 
\begin{defin}[{\cite[section 6]{GS3}}]
The arithmetic Grothendieck group $\ari{K}_{0}(X)$ associated to $X$ 
is the abelian group generated by $\widetilde{A}(X)$  and the isometry classes 
of hermitian bundles on $X$, with the following relations:
\begin{description}
\item[-\ ] $\widetilde{\ch}(\mtr{\cal E})=\mtr{E}\pr-\mtr{E}+\mtr{E}\prpr$ for every exact sequence 
$\mtr{\cal E}$ as above; 
\item[-\ ] $\eta=\eta\pr+\eta\prpr$ if $\eta\in\widetilde{A}(X)$ is the sum of two elements 
$\eta\pr$ and $\eta\prpr$.
\end{description}
\label{defarik}
\end{defin}
Notice that, by construction, there is an exact sequence of abelian groups
$$
\widetilde{A}(X)\ra\ari{K}_0(X)\ra K_0(X)\ra 0
$$
where the "forgetful" map $\ari{K}_0(X)\ra K_0(X)$ sends a hermitian 
bundle onto its underlying locally free sheaf and sends an element 
of $\widetilde{A}(X)$ to $0$. 

We shall now define a commutative ring structure on $\ari{K}_0(X)$. 
To this end, let us consider the group
 $\Gamma(X):=Z(X)\oplus\widetilde{A}(X)$. We equip it with 
the $\mN$-grading whose term of 
degree $p$ is $Z^{p,p}(X)\oplus\widetilde{A}
^{p-1,p-1}(X)$ if $p\geq 1$ and $Z^{0,0}(X)$ if $p=0$. 
We define an $\mR$-bilinear map $*$ from $\Gamma(X)\times\Gamma(X)$ to 
$\Gamma(X)$ via the formula 
$$
(\omega,\eta)*(\omega\pr,\eta\pr)=(\omega\wedge\omega\pr,\omega\wedge\eta\pr+
\eta\wedge\omega\pr+(\dbd\eta)\wedge\eta\pr).
$$
This map 
endows $\Gamma(X)$ with the structure of a commutative 
graded $\mR$-algebra (cf. \cite[Lemma 7.3.1, p. 233]{GS3}). 
Now let $\mtr{E}+\eta$ and $\mtr{E}\pr+\eta\pr$ 
be two generators of $\ari{K}_{0}(X)$; we define 
$$
(\mtr{E}+\eta)\otimes (\mtr{E}\pr+\eta\pr):=\mtr{E}\otimes\mtr{E}\pr+
[(\ch(\mtr{E}),\eta)*(\ch(\mtr{E}\pr),\eta\pr)].
$$
Here $[.]$ refers to the projection
on the second component of $\Gamma(X)$. 
Gillet and Soul\'e have shown in \cite[Th. 7.3.2]{GS3} that $\otimes$ is compatible with the defining relations of $\ari{K}_{0}(X)$ and defines 
a commutative ring structure on $\ari{K}_0(X)$.  

Now let $Y$ be another arithmetic variety.  Let $f:X\ra Y$ be any morphism. 
If $(E,h^E)+\eta$ is a generator of $\ari{K}_0(Y)$, we define 
$$f^*((E,h^E)+\eta)=(f^*(E),f^*_\mC h^E)+f^*_\mC(\eta),$$ where $f^*_\mC(\eta)$ is the pull-back 
of $\eta$ by $f_\mC$ as a differential form. It follows from 
the definitions that the just defined map $f^*$ 
descends to a morphism of commutative rings
$$
f^*:\ari{K}_0(Y)\ra\ari{K}_0(X).
$$
We shall call this morphism the {\it pull-back map} associated to $f$. 

We now turn to arithmetic Chow groups. 
We shall write $D^{p,p}(X)$ for the space of 
real currents of type $p,p$ on $X(\mC)$ on which 
$F_\infty^*$ acts by multiplication by $(-1)^p$.  
If $Z$ is a $p$-cycle on $X$, a Green current $g_Z$ for $Z$ is 
an element of $D^{p,p}(X)$ which satisfies the equation
$$
\dbd g_Z+\delta_{Z(\mC)}=\omega_Z
$$
where $\omega_Z$ is a differential form and $\delta_{Z(\mC)}$ is 
the Dirac current associated to $Z(\mC)$. 
\begin{defin}[{\cite[section 3]{GS2}}]
The arithmetic Chow group $\ari{\CH}^p(X)$ is 
the abelian group generated by the ordered pairs 
$(Z,g_Z)$, where $Z$ is a $p$-cycle on $X$ and $g_Z$ is a Green current for 
$Z(\mC)$, with the following relations: 
\begin{description}
\item[-\ ] $(Z,g_Z)+(Z',g_{Z'})=(Z+Z',g_{Z}+g_{Z'})$; 
\item[-\ ] $(\div(f),-\log|f|^2+\partial u+\mtr{\partial}v)=0$
\end{description}
where $f$ is a non-zero rational function defined on 
a closed integral subscheme of codimension $p-1$ in
$X$ and  
$u$ (resp. $v$) is a complex current of type $(p-2,p-1)$ (resp. 
$(p-1,p-2)$). 
\label{defarich}
\end{defin}
We shall write $\ari{\CH}^\cdot(X)$ for the direct sum 
$\oplus_{p\geqslant 0}\ari{\CH}^p(X)$. There is by construction 
a morphism of groups $\omega:\ari{CH}^\cdot(X)\ra Z(X)$, given 
by the formula $\omega((Z,g_Z)):=\omega_Z$. 
As for the arithmetic Grothendieck group, 
there is a natural exact sequence
$$
\widetilde{A}(X)\ra\ari{\CH}^\cdot(X)\ra\CH^\cdot(X)\ra 0
$$
where the "forgetful" map $\ari{\CH}^\cdot(X)\ra\CH^\cdot(X)$ sends 
a pair $(Z,g_Z)$ (as above) on $Z$. 

The maps $\widetilde{A}(X)\ra\ari{\CH}^\cdot(X)$ and 
$\widetilde{A}(X)\ra\ari{\CH}^\cdot(X)$ are usually 
both denoted by the letter $a$. To lighten formulae, we shall usually drop 
that letter in our computations. This is in the spirit of ordinary $K_0$-theory, 
where the brackets $[\cdot]$ (which map an object into the Grothendieck group) 
are often dropped in computations. 

The group 
$\ari{\CH}(X)_\mQ$ is equipped with a $\mZ$-bilinear pairing $\cdot$, such that 
$$
(Z,g_Z)\cdot(Z',g_{Z'})=(Z\cap Z',
g_Z\wedge\delta_{Z'(\mC)}+\omega_{Z}\wedge g_{Z'})
$$ 
if $Z,Z'$ are integral and meet properly in $X$;
the multiplicity of each component in $Z \cap Z'$ is 
given by Serre's $\rm{Tor}$ formula. 
 See \cite[Th. 4.2.3]{GS2} for the definition of the pairing in general. 
 It is proven in \cite{GS2} and \cite{Gubler} that this pairing makes the group $\ari{\CH}^\cdot(X)_\mQ$ into a commutative $\mN$-graded ring (the reference \cite{Gubler} fills a gap 
 in \cite{GS2}). 
Let now $f:X\ra Y$ be a morphism of 
arithmetic varieties. We can associate to $f$ a {\it pull-back map}
$$
f^*:\ari{\CH}^\cdot(Y)_\mQ\ra\ari{\CH}^\cdot(X)_\mQ.
$$
which is a morphism of $\mN$-graded rings. 
We shall describe this map under the hypothesis that 
$f$ is smooth over $\mQ$ and flat. 
Under this hypothesis, let $Z$ be a $p$-cycle on $Y$ and 
let $g_Z$ be a Green current for $Z$. Write 
$f^*Z$ for the pull-back of $Z$ to $X$ and $f^*g_Z$ for 
the pull-back of $g_Z$  to $X(\mC)$ as a current (which exists because $f_\mC$ is smooth). 
The rule which associates the pair 
$(f^*Z,f^*g_Z)$ to the pair $(Z,g_Z)$ descends to a morphism 
of abelian groups $\ari{\CH}^p(Y)\ra\ari{\CH}^p(X)$. 
The induced morphism $\ari{\CH}^p(Y)_\mQ\ra\ari{\CH}^p(X)_\mQ$ 
is the pull-back map. 
See \cite[sec. 4.4]{GS2}. 
 
There is a unique ring morphism 
$$
\ari{\ch}:\ari{K}_0(X)\ra\ari{\CH}^\cdot(X)_\mQ
$$
commuting with pull-back maps and such that
\begin{description} 
\item[(ch-1)] the formula $\ari{\ch}(\eta)=(0,\eta)$ holds, if $\eta\in\widetilde{A}(X)$; 
\item[(ch-2)] the formula $\ari{\ch}(\mtr{L})=\exp(\ari{\rm c}_1(\mtr{L}))$ holds, if $\mtr{L}=(L,h^L)$ is a hermitian line bundle on $X$; 
\item[(ch-3)] the formula $\omega(\ari{\ch}(\mtr{E}))=\ch(\mtr{E})$ holds for any hermitian vector bundle $\mtr{E}$ on $X$. 
\end{description}
Here the first Chern class $\ari{\rm c}_1(\mtr{L}) \in \ari{\CH}^1(X)$ of a hermitian line bundle
$\mtr{L}$ is defined as the class of $(\div\ (s),-\log\ h^L(s,s))$ for any choice of
a rational section $s$ of $L$ over $X$.
The fact 
that $-\log\ h^L(s,s)$ is a Green current for $\div\ (s)$ is implied 
by the Poincaré-Lelong formula (see \cite{Lelong}). 
The morphism $\ari{\ch}$ is called the {\it arithmetic Chern character} 
and is compatible with the traditional 
Chern character $K_0(X)\ra\CH^\cdot(X)_\mQ$ via the 
forgetful maps. See \cite[sec. 7.2]{GS3} for a proof of the existence 
and unicity of 
$\ari{\ch}$. 

The Todd class $\ari{\Td}(\mtr{E})$ of a hermitian vector bundle $\mtr{E}$
is defined similarly. It commutes with pull back maps and is multiplicative :
$$\ari{\Td}(\mtr{E}'\oplus\mtr{E}'')=\ari{\Td}(\mtr{E}')\ari{\Td}(\mtr{E}'') \, .$$
If $\mtr{L}$ is a hermitian line bundle, the formula
$$\ari{\Td}(\mtr{L}) = {\rm td}(\ari{\rm c}_1(\mtr{L}))$$
holds, where ${\rm td}(x)$ is the formal power series
${\rm td}(x) =  {xe^x\over e^x-1}\in\mQ[[x]]$.
The arithmetic 
Todd class is compatible with 
the usual Todd class via the forgetful maps. Furthermore, 
it has the following properties:
\begin{description}
\item[(Td-1)] if $\mtr{E}$ is a hermitian bundle on $X$, then 
$$
\omega(\ari{\Td}(\mtr{E}))=\Td(\mtr{E}); 
$$
\item[(Td-2)] if $\mtr{\cal E}$ is an exact sequence of hermitian 
bundles on $X$ as in Definition \ref{defarik}, then 
$$
\ari{\Td}(\mtr{E}'\oplus\mtr{E}'')-\ari{\Td}(\mtr{E})= \widetilde{\Td}(\mtr{\cal E}).
$$
\end{description}
For $\ari{\Td}$ to be uniquely defined, one still needs to give an expression
of the Todd class of the tensor product $ \mtr{E} \otimes \mtr{L}$
of a hermitian vector bundle with a hermitian line bundle.
See \cite[Th. 4.1, Th 4.8, par. 4.9]{GS3} for a proof. 

\section{The statement}

Recall that an {\it arithmetic variety} denotes a regular scheme, 
which is quasi-projective and flat  
over $\mZ$. 
Let $g:Y\ra B$ be a projective, flat  
morphism of arithmetic varieties, which is smooth 
over $\mQ$ (abbreviated p.f.s.r.). 

We shall first define a push-forward map $g_*:\ari{K}_0(Y)\ra\ari{K}_0(B)$. 
To this end, fix a conjugation invariant K\"ahler 
metric  $h_{Y}$ on $Y(\mC)$. Denote by 
$\omega_Y$ the corresponding Kähler form, given by the formula
$$
\omega_Y = i \sum_{\alpha , \beta} h_Y \left( \frac{\partial}{\partial z_{\alpha}} ,  \frac{\partial}{\partial z_{\beta}} \right) dz_{\alpha} \, d\bar z_{\beta}
$$
for any choice $(z_{\alpha})$ of local holomorphic coordinates.
Let $(E,h^E)$ be a hermitian bundle 
on $Y$, such that $E$ is $g$-acyclic. This means that $R^k g_*E=0$ if $k>0$ or 
equivalently that $H^k(Y_{\mtr{b}},E_{\mtr{b}})=0$ if $k>0$, for any 
geometric point $\mtr{b}\to B$. 
 The sheaf of modules $R^0 g_{*}E$ is then locally free by the semi-continuity 
 theorem. 
  Furthermore, in the holomorphic category, 
 the natural map  
 $$
 R^0 g_{\mC*}(E_\mC)_b\ra H^0(Y(\mC)_b,E(\mC)|_{Y(\mC)_b})
 $$
 is then an isomorphism for every point $b\in B(\mC)$. Here 
 $Y(\mC)_b$ denotes the (analytic) fiber of the morphism $g_\mC$ above 
 $b$. For every $b\in B(\mC)$, we endow $H^0(Y(\mC)_b,E(\mC)|_{Y(\mC)_b})$ with 
 the hermitian metric given by the formula
 $$
 \langle s,t\rangle_{L^2}:={1\over (2\pi)^{d_b}}\int_{Y(\mC)_b}h^E(s,t)\frac {\omega_Y^{d_b}}{d_b !}
 $$
 where $d_b:=\dim(Y(\mC)_b)$. It can be shown that these 
 metrics depend on $b$ in a $C^\infty$ manner (see \cite[p. 278]{BGV}) and thus define a hermitian metric on $(R^0 g_{*}E)_\mC$. We shall write $g_*h^E$ for this hermitian metric; 
 it is called the $L^2$-metric (obtained from $g_\mC$,  $h^E$ and $h_Y$).  
Apart from that, we shall write 
$T(h_{Y},h^{E})$ for the higher analytic torsion form determined by 
$(E,h^E)$, $g_\mC$ and $h_{Y}$.
The higher analytic torsion form is an element of 
$\widetilde{A}(B)$, which satisfies the equality
$$
\dbd T(h_Y,h^{E})=\ch((R^0 g_{*}E,g_{*}h^E))-\int_{Y(\mC)/B(\mC)}
\Td(\mtr{Tg}_\mC)\ch(\mtr{E}),
$$
where  
$\mtr{Tg}_{\mC}$ is the tangent bundle relatively to $g_\mC$, 
endowed with the hermitian metric induced by $h_Y$. 
For the definition of $T(h_{Y},h^{E})$ and for the proof of the last equality, we refer to 
\cite{BKo}. In \cite[Prop. 3.1]{R2}, it is shown that 
there is a unique group morphism 
$$
g_{*}:{\ari{K}_{0}}(Y)\ra {\ari{K}_{0}}(B)
$$ 
such that 
$$
g_{*}((E,h^E)+\eta)=(R^0 g_{*}E,g_{*}h^E)-T(h_{Y},h^{E})+\int_{Y(\mC)/B(\mC)}\Td(\mtr{Tg_\mC})\eta,
$$
where $\eta\in\widetilde{A}(Y)$ and $(E,h^E)$ is a hermitian bundle as above on $Y$. 
 We shall call the morphism $g_*$ 
  the {\it push-forward map} 
associated to $g$ and $h_{Y}$. 

There is also a {\it push-forward map} 
$$
g_*:\ari{\CH}^\cdot(Y)\ra\ari{\CH}^\cdot(B).
$$
This map is uniquely characterised by the fact that it is a group morphism and 
by the fact that 
 $$
 g_*((Z,g_Z))=(\deg(Z/g(Z))g(Z),\int_{Y(\mC)/B(\mC)}g_Z)
 $$ 
 for every
  integral closed subscheme $Z$ of $Y$ and for every Green current 
 $g_Z$ of $Z$. Here $\deg(Z/g(Z))$ is the degree of the corresponding 
 extension of function fields  if $\dim(g(Z))=\dim(Z)$ and 
 $\deg(Z/g(Z))=0$ otherwise. 
 
 Notice that 
 this push-forward map does not depend on the choice 
 of a Kähler metric on $Y(\mC)$, unlike the push-forward 
 map for arithmetic Grothendieck groups. 
 
Let now 
$$
\xymatrix{
Y \ar[rr]^i\ar[rd]^g & & P \ar[ld]^f\\
& B & \\
}
$$
be a factorisation of $g$ into a closed immersion $i$ and 
a projective smooth morphism $f$. Let $N$ be the normal 
bundle of the immersion $i$. Let 
$$
{\cal N}: 0\ra Tg_\mC\ra Tf_\mC\ra N_\mC\ra 0
$$
be the exact sequence associated to $i_\mC$. 
Endow as before $Tg_\mC$ with the metric induced by $h_Y$. 
Endow $Tf_\mC$ with some (not necessarily Kähler) hermitian metric extending the metric 
on $Tg_\mC$ and endow $N_\mC$ with the resulting quotient metric. These 
choices being made, we define
$$
\ari{\Td}(g)=\ari{\Td}(g,h_Y):=\ari{\Td}(i^*\mtr{Tf})\cdot\ari{\Td}^{-1}(\mtr{N})+
\widetilde{\Td}(\mtr{\cal N})\Td(\mtr{N})^{-1}\in\ari{\CH}^\cdot(Y)_\mQ.
$$
It is shown in \cite[Prop. 1, par. 2.6.2]{GS8} that 
the element $\ari{\Td}(g)$ depends only on $g$ and on the restriction 
of $h_Y$ to $Tg_\mC$. 

Before we state the Riemann-Roch theorem, we still have 
to define a characteristic class. 

\begin{defin}[{\cite[1.2.3, p. 25]{GS5}}]
The $R$-genus is the unique additive characteristic class 
defined
for a line bundle $L$ by the formula 
$$
R(L)=\sum_{m\ {\rm odd},\geq 1}(2\zeta\pr(-m)+\zeta(-m)(1+{1\over 2}+\dots
          +{1\over m}))c_{1}(L)^{m}/m!
$$
where $\zeta(s)$ is the Riemann zeta function.
\label{defR}
\end{defin}
In the definition \ref{defR}, it is understood as usual that $R$ is defined 
for any $C^\infty$-vector bundle on a $C^\infty$-manifold and that it has 
values in ordinary (de Rham) cohomology with complex coefficients. 
This being said, let 
$$
H_\red^\cdot(Y):=\sum_{p\geqslant 0}Z^{p,p}(Y)/(A^{p,p}(Y)\cap\Im\ d).
$$
By construction, there is an inclusion 
$H_\red^\cdot(Y)\subset H^\cdot(Y(\mC),\mC)$. If $E$ is a vector 
bundle on $Y$, then $R(E_\mC)$ can be computed via the formulae 
of Chern-Weil using a connection of type $(1,0)$. The local 
curvature matrices associated to such connections are of 
type $(1,1)$; this shows that $R(E_\mC)\in H_\red^\cdot(Y)$. 
On the other hand, also by construction,
 there is a natural map $H_\red^\cdot(Y)\ra\widetilde{A}(Y)$. 
 Hence we may (and shall) consider that $R(E_\mC)\in\widetilde{A}(Y)$. 
 Similar remarks apply to any other characteristic class. 
\begin{theor}[arithmetic Riemann-Roch theorem]
Let $y\in\ari{K}_{0}(Y)$. The equality 
$$
\ari{\ch}(g_{*}(y))=g_{*}(\ari{\Td}(g)\cdot (1-a(R(Tg_\mC )))\cdot\ari{\ch}(y))
$$
holds in $\ari{\CH}^\cdot(B)_\mQ$.
\label{arr}
\end{theor}
In \cite{GS8}, it is proved that the equality in Theorem \ref{arr} holds 
after projection of both sides of the equality on 
$\ari{\CH}^1(B)_\mQ$. 

\section{The Proof}

In this section, we shall prove Theorem \ref{arr}. 
The structure of the proof is as follows. 
In the first subsection, we prove various properties 
of the (putatively non vanishing) difference between the two sides of 
the asserted equality. Let us call this difference the error term. 
We first prove that the error term 
is independent of all the involved hermitian metrics (Lemmata 
\ref{deltadiff}, \ref{indepmet} and \ref{indepkahl}), using 
Bismut-Koehler's anomaly formulae for the analytic torsion form. 
We then proceed to prove 
that it is invariant under immersions (Theorem \ref{imerr}). The proof of this 
fact relies on two difficult results, which are proved elsewhere: 
the arithmetic Riemann-Roch for closed immersions (Theorem \ref{rrim}), which 
is a generalisation of Arakelov's adjunction formula and 
Bismut's immersion formula (Theorem \ref{BLA}). This last result is the most 
difficult part of the proof of the arithmetic Riemann-Roch 
theorem and is of a purely analytic nature. In the 
third section, we show that the error term vanishes in 
the special case of relative 
projective spaces; this is is shown to be either the consequence of 
the article \cite{GS8}, where the arithmetic Riemann-Roch theorem is proved in degree $1$  
or of the more recent article by Bost \cite{JBB1}, where explicit resolutions of 
the diagonal are used. Finally, in the third subsection, we show 
that the error term always vanishes. This is achieved 
 by reduction to the case of 
relative projective spaces, using the invariance of the error term 
under immersions.

\subsection{Properties of the error term}

Before beginning with the study of the properties of the error 
term of the arithmetic Riemann-Roch theorem, we shall recall 
a few results  on direct images in arithmetic Chow and $K_0$-theory. 

Let $i:Y\hookrightarrow P$ be a closed immersion of 
arithmetic varieties. Let $\eta$ be a locally free sheaf on $Y$ and let 
$$
\Xi:\ 0\to \xi_{m}\to\xi_{m-1}\to\dots \xi_{0}\to i_{*}\eta\to 0
$$
be a resolution of $i_*\eta$ by locally free sheaves on $P$. Denote by  
$N$ be the normal bundle of the immersion $i$. 
 Let 
$F:=\oplus_{l=0}^{m}H^{l}(i^*\Xi)$. 
 There is a canonical isomorphism of graded bundles 
$F\simeq \oplus_{l=0}^{rk(N)}\Lambda^{l}(N^{\vee})\otimes 
\eta$ (see for instance \cite[Lemme 2.4 and Prop. 2.5, i'), exposé VII]{SGA6}). 
Both of the latter graded bundles 
carry natural metrics, if $N$, $\eta$ and the 
$\xi_i$ are endowed with metrics. 
We shall say that 
hermitian metrics on the bundles $\xi_{i}$ satisfy Bismut's assumption (A) 
with respect to the hermitian metrics on $N$ and $\eta$ if 
the isomorphism $i^{*}F\simeq \oplus_{i=0}^{rk(N)}\Lambda^{i}(N^{\vee})\otimes 
\eta$ also identifies the metrics. 
It is proved in \cite{BGS2} that if metrics on 
$N$ and $\eta$ are given, there always exist metrics 
on the $\xi_{i}$ such that this assumption is satisfied.   
We now equip $N$ and $\eta$ with arbitrary hermitian 
metrics and we suppose that the $\xi_{i}$ are endowed 
with hermitian metrics such that Bismut's condition (A) is satisfied with respect 
to the metric on $N$ and $\eta$.  
The {\it singular Bott-Chern current} of  $\mtr{\Xi}$ is 
an element $T(h^{\xi_{\cdot}})$ of $\oplus_{p\geqslant 0} D^{p,p}(X)$ 
satisfying the equation
$$
\dbd T(h^{\xi_{\cdot}})=i_{*}(\Td^{-1}(\mtr{N})\ch(\mtr{\eta}))-
\sum_{i=0}^{m}(-1)^{i}\ch(\mtr{\xi}_{i})
$$
(see 
\cite[Th. 2.5, p. 266]{BGS2}).
 Here $i_{*}$ refers to the pushforward of currents.

We now suppose given a commutative 
diagram 
$$
\xymatrix{
Y \ar[rr]^i\ar[rd]^g & & P \ar[ld]^f\\
& B & \\
}
$$
where $g$ is a  p.f.s.r. morphism and $f$ is projective and smooth.
Endow 
$P(\mC)$ with a Kähler metric $h_P$ and $Y(\mC)$ with the restricted metric $h_Y$. 
As before Definition \ref{defR}, consider the sequence 
$$
{\cal N}: 0\ra Tg_\mC\ra Tf_\mC\ra N_\mC\ra 0.
$$
Endow $Tf_\mC$ (resp. $Tg_\mC$) with the metric induced from $h_P$ 
(resp. $h_Y$). With these conventions, we shall suppose from now on that the metric on 
$N_\mC$  is the quotient metric induced 
from the map $Tf_\mC\ra N_\mC$ in the sequence $\cal N$. 
 The following 
result is proved in \cite[Th. 4.13]{BGS3}. 
\begin{theor}[arithmetic Riemann-Roch theorem for closed immersions]
Let $\alpha\in\ari{CH}^\cdot(P)$. The current $$\int_{P(\mC)/B(\mC)}
\omega(\alpha)T(h^{\xi})$$ is then a differential form and the equality 
$$
f_{*}(\alpha\cdot\ari{\ch}(\xi_{\cdot}))=g_{*}(i^{*}(\alpha)\cdot
\ari{\Td}^{-1}(\mtr{N})\cdot\ari{\ch}(\mtr{\eta}))-\int_{P(\mC)/B(\mC)}
\omega(\alpha)T(h^{\xi})
$$
is satisfied in $\ari{\CH}^\cdot(B)_\mQ$. 
\label{rrim}
\end{theor}
Notice that if one applies the forgetful map to both 
sides of the last equality, one obtains a consequence of 
the Grothendieck-Riemann-Roch  
theorem for closed immersions. 

Suppose from now on that the $\xi_i$ are $f$-acyclic and that 
$\eta$ is $g$-acyclic. 
The next theorem is a $\ari{K}_0$-theoretic translation of a difficult 
result of Bismut, often called Bismut's immersion theorem. The translation 
is made in \cite[Th. 6.6]{R2}. Bismut's immersion theorem 
is proved in \cite{SFP}.
\begin{theor}
The equality
\begin{eqnarray*}
&&g_{*}(\mtr{\eta})-\sum_{i=0}^{m}(-1)^{i}f_{*}(\mtr{\xi}_{i})=\\
&=&
\int_{Y(\mC)/B(\mC)}\ch(\eta_\mC)R(N_\mC)\Td(Tg_\mC)+\int_{P(\mC)/B(\mC)}T(h^{\xi_{\cdot}})\Td(\mtr{Tf})\\
&+&\int_{Y(\mC)/B(\mC)}\ch(\mtr{\eta})\widetilde{\Td}(\mtr{\cal N})\Td^{-1}(\mtr{N})
\end{eqnarray*}
holds in $\ari{K}_{0}(B)$.
\label{BLA}
\end{theor}
Notice that the complex $R^0 f_*(\Xi)$ is exact with our hypotheses. 
Hence we have $g_{*}({\eta})-\sum_{i=0}^{m}(-1)^{i}f_{*}({\xi}_{i})=0$ 
in $K_0(B)$. This is the equality to which the Theorem \ref{BLA} 
reduces after application of the forgetful maps. 
We shall also need 
the following theorem, which studies the dependence of the analytic torsion form 
on $h_Y$:
\begin{theor}
Let $h'_Y$ be another Kähler metric on $Y(\mC)$.  Let $h^{'Tg}$ be  the metric induced 
on $Tg_\mC$ by $h'_Y$. The 
identity 
$$
T(h'_Y,h^{\eta})-T(h_Y,h^{\eta})=\widetilde{\ch}(g^{h_Y}_{*}h^{\eta},g^{h'_Y}_{*}h^{\eta})-\int_{Y(\mC)/B(\mC)}\widetilde{\Td}(h^{Tg},h^{'Tg})\ch(\mtr{\eta}).
$$
holds in $\widetilde{A}(B)$.
\label{ArComp}
\end{theor}
Here $\widetilde{\Td}(h^{Tg},h^{'Tg})$ refers to the Todd secondary class of the sequence
$$0\ra 0\ra Tg_\mC\ra Tg_\mC\ra 0, $$ where the first non-zero term is endowed with the metric 
$h^{Tg}$ and the second non-zero term with the metric $h^{'Tg}$. The term $\widetilde{\ch}(g^{h_Y}_{*}h^{\eta},g^{h'_Y}_{*}h^{\eta})$ 
is the Bott-Chern secondary class of the  sequence $$0\ra 0\ra R^0 g_{*}\eta\ra R^0 g_{*}\eta\ra 0, $$ 
where the first non-zero term carries the metric obtain by integration along the fibers 
with the volume form coming from $h_Y$ and the second non-zero term the metric obtain by integration along the fibers 
with the volume form coming from $h'_Y$. 
For the proof, we refer to \cite[Th. 3.10, p. 670]{BKo}. 

We are now ready to study the {\it error term}
$$
\delta(y,g,h_{Y}):=\ari{\ch}(g_{*}(y))-g_{*}(\ari{\Td}(g)\cdot(1-R(Tg_\mC ))\cdot\ari{\ch}(y))
$$
of the arithmetic Riemann-Roch theorem. Notice that by construction 
$$\delta(y'+y'',g,h_Y)=\delta(y',g,h_Y)+\delta(y'',g,h_Y)$$ for all 
$y',y''\in\ari{K}_0(Y)$. 
\begin{lemma} $\delta(y,g,h_Y)=0$  if $y$ is represented by a differential 
form.
\label{deltadiff}
\end{lemma}
\beginProof
This follows directly from the definitions.
\qed\par
\begin{lemma}
Let $\mtr{E},\mtr{E}'$ be hermitian vector bundles on $Y$ such that 
$E\simeq E'$. Then we have $\delta(\mtr{E},g, h_Y)=\delta(\mtr{E}\pr,g,h_Y)$. 
\label{indepmet}
\end{lemma}
\beginProof
We have $\delta(\mtr{E},g)-\delta(\mtr{E}\pr,g)=\delta(\mtr{E}-\mtr{E}\pr,g)$ 
and from the definition of arithmetic $K_0$-theory, the element 
$\mtr{E}-\mtr{E}\pr$ is represented by a differential form. Hence we can apply the last 
lemma. 
\endProof
\begin{lemma}
Let $h'_Y$ be another K\"ahler metric on $Y$. Then 
$\delta(y,g,h_{Y})=\delta(y,g,h_{Y}\pr)$.
\label{indepkahl}
\end{lemma}
\beginProof
We shall use Theorem \ref{ArComp}. Using Lemma \ref{deltadiff} and 
the fact that every locally free  sheaf on $Y$ has a finite resolution by 
$g$-acyclic locally free sheaves, we see that  
we may assume without loss of generality that 
$y=\mtr{E}=(E,h^E)$, where $\mtr{E}$ a hermitian bundle on $Y$ such that 
$E$ is $g$-acyclic. We now compute 
\begin{eqnarray*}
&&\delta(\mtr{E},g,h_{Y})-\delta(\mtr{E},g,h'_{Y})\stackrel{(1)}{=}\\
&=&\ari{\ch}((g_{*}E,g_{*}^{h_{Y}}h^{E}))-T(h_Y,h^{E})-
g_{*}(\ari{\Td}(g,h_Y)\cdot(1-R(Tg_\mC ))\ari{\ch}(\mtr{E}))-\\
&&\ari{\ch}((g_{*}E,g_{*}^{h'_Y}h^E))+T(h'_Y,h^{E})+
g_{*}(\ari{\Td}(g,h'_Y)\cdot(1-R(Tg_\mC ))\ari{\ch}(\mtr{E}))\stackrel{(2)}{=}\\
&=&-\widetilde{\ch}(g_{*}^{h_Y}h^{E},g_{*}^{h'_Y}h^{E})+
T(h'_Y,h^{E})-
T(h_Y,h^{E})\\
&&+g_{*}(\widetilde{\Td}(h^{Tg},h^{'Tg})\cdot(1-R(Tg_\mC ))
\ari{\ch}(\mtr{E}))\\
&\stackrel{(3)}{=}&
T(h'_Y,h^{E})-T(h_Y,h^{E})\\
&&-\Big(
\widetilde{\ch}(g_{*}^{h_Y}h^{E},
g_{*}^{h'_Y}h^{E})-\int_{Y(\mC)/B(\mC)}\ch(\mtr{E})\widetilde{\Td}(h^{Tg},h^{'Tg})\Big)
\stackrel{(4)}{=}0 \, .
\end{eqnarray*}
The equality (1) follows from the definitions. The equality (2) 
is justified by the property (ch-1) of the arithmetic Chern character 
and by  \cite[par. 2.6.2, Prop.1,(ii)]{GS8}, which 
implies that 
$$
\ari{\Td}(g,h'_Y)-\ari{\Td}(g,h_Y)=\widetilde{\Td}(h^{Tg},h^{'Tg}) 
$$
(notice that this follows from (Td-2) if $g$ is smooth). 
From the definition of the ring structure 
of $\ari{\CH}^\cdot(X)_\mQ$ (see after Definition \ref{defarich}), 
we obtain (3). The equality (4) is the content of Theorem \ref{ArComp}.
\endProof
In view of the Lemma \ref{indepkahl}, we shall from now on  drop the 
reference to the Kähler metric and 
write $\delta(y,g)$ for the error term. Notice that 
the Lemmata \ref{deltadiff} and \ref{indepmet} imply that 
$\delta(y,g)$ depends only on the image of $y$ in 
$K_0(Y)$. This justifies writing $\delta(E,g)$ for 
$\delta(\mtr{E},g)$ if $\mtr{E}$ is a hermitian bundle on $Y$. 

The following theorem studies the compatibility 
of the error term with the immersion $i$ and is the core 
of the proof of the arithmetic Riemann-Roch theorem. 
\begin{theor}
The equality  
$$
\sum_{i=0}^{m}(-1)^{i}\delta({\xi}_{i},f)=\delta({\eta},g)
$$
holds. 
\label{imerr}
\end{theor}
\beginProof 
 Using Theorem \ref{BLA}, we compute 
\begin{eqnarray*}
\sum_{i=0}^{m}(-1)^{i}\delta({\xi}_{i},f)&=&\sum_{i=0}^{m}
(-1)^{i}(\ari{\ch}(f_{*}(\mtr{\xi}_{i}))-
f_{*}(\ari{\Td}(\mtr{Tf})\cdot(1-R(Tf_\mC ))\cdot\ari{\ch}(\mtr{\xi}_{i})))=\\
&=&
\ari{\ch}(g_{*}(\mtr{\eta}))-\int_{Y(\mC)/B(\mC)}\ch(\eta_\mC)R(N_\mC)\Td(Tg_\mC)-\int_{P(\mC)/B(\mC)}T(h^{\xi_{\cdot}})\Td(\mtr{Tf})\\
&&
-\int_{Y(\mC)/B(\mC)}\ch(\mtr{\eta})\widetilde{\Td}(\mtr{\cal N})\Td^{-1}(\mtr{N})\\
&&
-\sum_{i=0}^{m}(-1)^{i}f_{*}(\ari{\Td}(\mtr{Tf})\cdot\ari{\ch}(\mtr{\xi}_{i}))+
\sum_{i=0}^{m}(-1)^{i}\int_{P(\mC)/B(\mC)}\Td(Tf_\mC)
\ch({\xi}_{i,\mC})R(Tf_\mC).
\end{eqnarray*}
Now by the definition of the arithmetic tangent element, we have 
$$
\widetilde{\Td}(\mtr{\cal N})\Td^{-1}(\mtr{N})+\ari{\Td}(i^*\mtr{Tf})\cdot\ari{\Td}^{-1}(\mtr{N})=\ari{\Td}(g)
$$
and hence 
$$
\int_{Y(\mC)/B(\mC)}\ch(\mtr{\eta})\widetilde{\Td}(\mtr{\cal N})\Td^{-1}(\mtr{N})=g_*(\ari{\Td}(g)\cdot\ari{\ch}(\mtr{\eta}))-g_*(\ari{\ch}(\mtr{\eta})\cdot
\ari{\Td}(i^*\mtr{Tf})\cdot\ari{\Td}^{-1}(\mtr{N})).
$$
Hence, using Theorem \ref{rrim} with $\alpha = \ari{\Td}(\mtr{Tf})$
 and property (Td-1) of the arithmetic Todd 
class, we obtain that 
\begin{eqnarray*}
\sum_{i=0}^{m}(-1)^{i}\delta({\xi}_{i},f)
&=&
\ari{\ch}(g_{*}(\mtr{\eta}))-g_*(\ari{\Td}(g)\cdot\ari{\ch}(\mtr{\eta}))-
\int_{Y(\mC)/B(\mC)}\ch(\eta_\mC)R(N_\mC)\Td(Tg_\mC)+\\
&&
\sum_{i=0}^{m}(-1)^{i}\int_{P(\mC)/B(\mC)}\Td(Tf_\mC)
\ch({\xi}_{i,\mC})R(Tf_\mC).
\end{eqnarray*}
Furthermore, using the Grothendieck-Riemann-Roch theorem with values 
in singular cohomology, we compute that 
 $$
\sum_{i=0}^m(-1)^{i}\ch(\xi_{i,\mC})=i_{*}(\Td^{-1}(N_\mC)\ch(\eta_\mC))
$$
where $i_*$ is the direct image in singular cohomology. Hence, using 
the multiplicativity of the Todd class and the additivity of the 
$R$-genus, 
\begin{eqnarray*}
&&\sum_{i=0}^{m}(-1)^{i}\int_{P(\mC)/B(\mC)}\Td(Tf_\mC)
\ch({\xi}_{i,\mC})R(Tf_\mC)-\int_{Y(\mC)/B(\mC)}\ch(\eta_\mC)R(N_\mC)\Td(Tg_\mC)\\
&=&\int_{Y(\mC)/B(\mC)}\ch(\eta_\mC)\Td(Tg_\mC)(R(Tf_\mC)-R(N_\mC))=\int_{Y(\mC)/B(\mC)}\ch(\eta_\mC)\Td(Tg_\mC)R(Tg_\mC).
\end{eqnarray*}
Thus
$$
\sum_{i=0}^{m}(-1)^{i}\delta({\xi}_{i},f)=\ari{\ch}(g_{*}(\mtr{\eta}))-g_*(\ari{\Td}(g)\cdot(1-R(Tg_\mC))\cdot\ari{\ch}(\mtr{\eta}))=\delta({\eta},g)
$$
which was the claim to be proved.
\qed\par

\subsection{The case of relative projective spaces}

Suppose now that $B$ is an arithmetic variety and that 
$Y=\mP^r_B\simeq\mP^r_\mZ\times B$ is some relative 
projective space over $B$ ($r\geqslant 0$).
 Let $g:Y\ra B$ (resp. $p:Y\ra\mP^r_\mZ$) 
be  the natural 
projection.  Endow $Y(\mC)$ with the product $h_Y$ of 
the standard Fubini-Study Kähler metric on $\mP^r(\mC)$ 
with a fixed (conjugation invariant) Kähler metric on $B(\mC)$. Let 
$\omega_{\mP^r}$ and $\omega_B$ be the corresponding Kähler 
forms on $\mP^r(\mC)$ and $B(\mC)$. 
Let $k\in\mZ$ and let $\mtr{\cal O}(k)$ be the $k$-th tensor power of the tautological bundle on $\mP_\mZ^r$, endowed with the Fubini-Study metric. 
We shall write $\tau(\mtr{\cal O}(k))$ for the analytic torsion form 
of $\mtr{\cal O}(k)$ with respect to the map from 
$\mP(\mC)$ to the point. The form $\tau(\mtr{\cal O}(k))$ is 
in this case a real number (which coincides with the Ray-Singer torsion 
of $\mtr{\cal O}(k)_\mC$).  We shall need the 
\begin{lemma}
Let $\mtr{E}$ be a hermitian vector bundle 
on $B$. 
Let $\mtr{V}:=g^*(\mtr{E})\otimes p^*(\mtr{\cal O}(k))$. Then 
the equality $$T(h_Y,h^V)=\ch(\mtr{E})\tau(\mtr{\cal O}(k))$$ 
holds for the analytic torsion form $T(h_Y,h^V)$ of $\mtr{V}$ with respect
 to $g$ and $h_Y$. 
 \label{projtor}
\end{lemma}
\beginProof
See \cite[Lemma 7.15]{R2}. In that reference, it is assumed that $k>>0$ but this assumption 
is not used in the proof and is thus not necessary.
\endProof
We shall also need the following projection formulae.
\begin{prop}
Let $b\in\ari{K}_0(B)$ and $a\in\ari{K}_0(Y)$. Then the projection 
formula 
$$
g_*(a\otimes g^*(b))=g_*(a)\otimes b
$$
holds in $\ari{K}_0(B)$.

Similarly, let $b\in\ari{\CH}^\cdot(B)_\mQ$ and 
$a\in\ari{\CH}^\cdot(Y)_\mQ$. Then the projection formula
$$
g_*(a\cdot g^*(b))=g_*(a)\cdot b
$$
holds in $\ari{\CH}^\cdot(B)_\mQ$.
\label{projform}
\end{prop}
\beginProof
For the first formula, see \cite[Prop. 7.16]{R2}. For the second one, see 
\cite[Th. in par. 4.4.3]{GS2}. 
\endProof
\begin{prop} 
We have $\delta({\cal O}_Y,g)^{[0]}=\delta({\cal O}_Y,g)^{[1]}=0$.
\end{prop}
\beginProof
To prove that $\delta({\cal O}_Y,g)^{[0]}=0$, notice that 
the forgetful map $\ari{\CH}^0(B)\ra\CH^0(B)$ is an isomorphism 
by construction. Hence the assertion that $\delta({\cal O}_Y,g)^{[0]}=0$ 
follows from the Grothendieck-Riemann-Roch theorem. 
The fact that $\delta({\cal O}_Y,g)^{[1]}=0$ is a special case of 
\cite[Th. 7]{GS8}. A different proof is given in \cite[par. 4.2]{JBB1}. 
\endProof
\begin{cor}
We have $\delta({\cal O}_Y,g)=0$. 
\label{trivbund}
\end{cor}
\beginProof 
Endow ${\cal O}_Y$ with the trivial metric. 
Recall that
$$
\delta({\cal O}_Y,g)=\ari{\ch}(g_*(\mtr{\cal O}_Y))
-g_*(\ari{\Td}(g)(1-R(Tg_\mC))).
$$
We shall first show that $\ari{\ch}(g_*(\mtr{\cal O}_Y))$ has 
no components of degree $>1$ in $\ari{\CH}^\cdot(B)_\mQ$. 

Notice that for every $l>0$, we have $R^l g_*({\cal O}_Y)=0$ and that  
the isomorphism $R^0 g_*({\cal O}_Y)\simeq{\cal O}_B$ given 
by adjunction is an isomorphism. 
We shall compute the $L^2$-norm of the section $1$, which trivialises 
$R^0 g_*({\cal O}_Y)$. For $b\in B(\mC)$, we compute  
\begin{eqnarray*}
\langle1_b,1_b\rangle_{L^2}={1\over (2\pi)^r r!}\int_{\mP^r(\mC)}
i_b^*(g^*(\omega_B)+p^*(\omega_{\mP^r}))^r = {1\over (2\pi)^r r!}\int_{\mP^r(\mC)}
\omega_{\mP^r}^r
\end{eqnarray*}
where $i_b:\mP^r(\mC)\hookrightarrow Y$ is the embedding 
of the fiber of the map $g_\mC$ above $b$. The first equality 
holds by definition. The second one is justified by the binomial 
formula and by the fact that $i_b^*\cdot g^*(\omega_B)=0$, since the image of 
$g\circ i_b$ is the point $b$. 
We thus see that $\langle 1_b,1_b\rangle_{L^2}$ is independent of 
$b\in B(\mC)$. This shows that $R^0 g_*({\cal O}_Y)$ 
endowed with its $L^2$-metric is the trivial bundle endowed with 
a constant metric. Aside from that, by Lemma \ref{projtor} we have 
$T(h_Y,h^{{\cal O}_Y})=\tau(\mtr{\cal O}_{\mP^r(\mC)})$. 
Hence the differential form $T(h_Y,h^{{\cal O}_Y})$ is a constant 
function on $B(\mC)$. Now, using the definition of 
the push-forward map in arithmetic $K_0$-theory (see Section 3), we 
compute that 
$$
\ari{\ch}(g_*(\mtr{\cal  O}_Y))=\ari{\ch}((R^0 g_*({\cal O}_Y),g_*h^{{\cal O}_Y})
-\tau(\mtr{\cal O}_{\mP^r(\mC)}). 
$$
Using the property (ch-2) of the Chern character, we see that 
$\ari{\ch}((R^0 g_*({\cal O}_Y),g_*h^{{\cal O}_Y}))$ 
has no components of degree $>1$. 
We can thus conclude that $\ari{\ch}(g_*(\mtr{\cal  O}_Y))$ 
has no components of degree $>1$. 

We shall now show that  $g_*(\ari{\Td}(g)(1-R(Tg_\mC)))$ has no 
components of degree $>1$. This will conclude the proof. 
By construction, we have 
$$
g_*(\ari{\Td}(g)\cdot(1-R(Tg_\mC)))=g_*\big(p^*(\ari{\Td}(\mtr{T\mP^r_\mZ})\cdot(1-R(T\mP^r_\mC)))\big).
$$
Now, since $\mP^r_\mZ$ has dimension $r+1$, the element 
$p^*(\ari{\Td}(\mtr{T\mP^r_\mZ})(1-R(T\mP^r_\mC)))$ 
has no component of degree $>r+1$. Hence the element 
$$
g_*\big(p^*(\ari{\Td}(\mtr{T\mP^r_\mZ})\cdot(1-R(T\mP^r_\mC)))\big)
$$
has no component of degree $>1$. 

We thus see that $\delta({\cal O}_Y,g)$ has no 
components of degree $>1$. We can now conclude 
the proof using the last Proposition.
\endProof
\begin{cor}
We have $\delta(y,g)=0$ for all $y\in\ari{K}_0(Y)$.
\label{arrproj}
\end{cor}
\beginProof
We know that $K_0(Y)$ is generated by elements of the form 
$g^*({E})\otimes p^*({\cal O}(k))$, where $E$ is a vector bundle 
on $B$ and $k\geqslant 0$ (see \cite[exp. VI]{SGA6}). In view of 
Lemma \ref{deltadiff} and Lemma \ref{indepmet}, we are thus reduced to prove 
that $\delta(g^*({E})\otimes p^*({\cal O}(k)),g)=0$. 
Now in view of Proposition \ref{projform} and the fact that 
the arithmetic Chern character  is multiplicative and commutes 
with pull-backs, we have 
$$
\delta(g^*({E})\otimes p^*({\cal O}(k)),g)=\delta(p^*({\cal O}(k)),g)\cdot\ari{\ch}((E,h^E)).
$$
for any hermitian metric $h^E$ on $E$. 
We are thus reduced to prove that 
$\delta(p^*({\cal O}(k)),g)=0$ (for $k\geqslant 0$). 
In order to emphasize the fact that $g$ depends only on $B$ and $r$,
 we shall write $\delta(p^*({\cal O}(k)),B,r)$ for $\delta(p^*({\cal O}(k)),g)$ 
 until the end of the proof. 
Now note that $\delta(p^*({\cal O}(k)),B,r)=0$ if $r=0$. 
Furthermore, $\delta(p^*({\cal O}(k)),B,r)=0$ if $k=0$ by Corollary \ref{trivbund}. 
By induction, we may thus assume that $k,r>0$ and that 
$\delta(p^*({\cal O}(k')),B,r')=0$ for all $k',r'\in\mN$ such that 
$k^{'2}+r^{'2}<k^2+r^2$. 
Now recall that there is an exact sequence of coherent sheaves  
$$
0\ra{\cal O}(-1)\ra {\cal O}_{\mP^r_\mZ}\ra j_*{\cal O}_{\mP^{r-1}_\mZ}\ra 0
$$
where $j$ is the immersion of $\mP^{r-1}_\mZ$ into 
$\mP^r_\mZ$ as the hyperplane at $\infty$. If we tensor 
this sequence with ${\cal O}(k)$, we obtain the sequence
$$
0\ra{\cal O}(k-1)\ra {\cal O}(k)\ra j_*{\cal O}_{\mP^{r-1}_\mZ}(k)\ra 0.
$$
If we apply Theorem \ref{imerr} to this sequence, we see that 
the equalities\linebreak$\delta(p^*({\cal O}(k-1)),B,r)=0$ 
 and $\delta(p^*({\cal O}(k)),B,r-1)=0$ together imply the 
 equality $\delta(p^*({\cal O}(k)),B,r)=0$. The first two equalities 
 hold by induction, so this concludes the proof. 
\endProof

\subsection{The general case}

To conclude the proof of the arithmetic Riemann-Roch theorem, 
we consider again the case of a general p.f.s.r morphism of 
arithmetic varieties $g:Y\ra B$. We want to prove 
that $\delta(y,g)=0$ for every  $y\in{K}_0(Y)$. 
Since $K_0(Y)$ is generated by 
$g$-acyclic bundles, we may assume that 
$y=\eta$, where $\eta$ is a $g$-acyclic vector bundle. 

Now notice that, by assumption, there is 
an $r\in\mN$ and a commutative diagram 
$$
\xymatrix{
Y \ar[rr]^i\ar[rd]^g & & \mP^r_B \ar[ld]_f\\
& B & \\
}
$$
where $f$ is the natural projection and $i$ is a closed immersion. 
Choose a resolution 
$$
0\to \xi_{m}\to\xi_{m-1}\to\dots \xi_{0}\to i_{*}\eta\to 0
$$
of $i_*\eta$ by $f$-acyclic locally free sheaves $\xi_i$ on $\mP^r_B$. 
Corollary \ref{arrproj} implies that $\delta(\xi_i,f)=0$. From this and 
Theorem \ref{imerr} we deduce that 
$\delta(\eta,g)=0$ and this concludes the proof of the arithmetic 
Riemann-Roch theorem.

\bibliographystyle{plain}

 \end{document}